\documentclass[twoside]{article}
\title{Pr\"ufer-Like Conditions in Subring Retracts and Applications}
\date{}
\usepackage{amsfonts}
\usepackage{amsmath}
\usepackage{amssymb}
\usepackage{latexsym}

\newtheorem{thm}{\bf Theorem}[section]

\newtheorem{prop}[thm]{\bf Proposition}

\newtheorem{exmp}[thm]{\bf Example}

\catcode`\ç=13
\defç{\c{c}}
\catcode`\é=13
\defé{\'e}
\catcode`\à=13
\defà{\`a}
\catcode`\è=13
\defè{\`e}
\catcode`\â=13
\defâ{\^a}
\catcode`\ù=13
\defù{\`u}
\catcode`\ê=13
\defê{\^e}
\catcode`\î=13
\defî{\^\i}
\catcode`\ô=13
\defô{\^o}
\newcommand{\field}[1]{\mathbb{#1}}

\newcommand{\Z }{\field{Z}}

\def\proof{{\parindent0pt {\bf Proof.\ }}}

\newcommand{\cqfd}
{\hspace{1cm}
\rule{2mm}{2mm}%
\medbreak%
\par%
}
\begin{document}
\thispagestyle{empty}

\maketitle \vspace*{-1.5cm}

\begin{center}
{\large\bf Chahrazade Bakkari $^1$, Najib Mahdou $^1$,  and Hakima Mouanis $^2$}

\bigskip

$^1$ Department of Mathematics, Faculty of Science and Technology of Fez,\\ 
Box 2202, University S. M. Ben Abdellah Fez, Morocco,\\
cbakkari@hotmail.com \\ mahdou@hotmail.com\\
 $^{2}$  Department of Mathematics, Faculty of Science of Rabat,\\ 
University Mohammed V Rabat, Morocco,\\
hmouanis@yahoo.fr\\

\end{center}

\bigskip\bigskip

\noindent{\large\bf Abstract.}  In this paper, we consider five
possible extensions of the Pr\"ufer domain notion to the case of
commutative rings with zero divisors. We investigate the 
transfer of these Pr\"ufer-like properties between a commutative ring and its 
 subring retract. Our results generate new
families of examples of rings subject to a given Pr\"ufer-like
conditions.\bigskip

\small{\noindent{\bf Key Words.} Pr\"ufer rings, Gaussian rings,
arithmetical rings, weak global dimension of rings, semihereditary
rings, subring retract, trivial ring extensions, Nagata rings.
\bigskip\bigskip

\begin{section}{Introduction}   Throughout this paper all rings are commutative with identity
element and all modules are unital. \\

In his article [25], Pr\"ufer introduced a new class of integral
domains, namely those domains $R$ in which all finitely generated
ideals are invertible. Through the years, Pr\"ufer domains
acquired a great many equivalent characterizations, each of wich
can, and was, extended to rings with zero-divisors in a number of
ways. More precisely, we consider
the following Pr\"ufer-like properties on a commutative ring ([3] and [4]):  \\

\indent {\bf (1)} $R$ is semihereditary, i.e., every finitely generated ideal is projective. \\
\indent {\bf (2)} The weak global dimension of $R$ is at most one. \\
\indent {\bf (3)} $R$ is an arithmetical ring, i.e., every
finitely generated ideal
is locally principal. \\
\indent {\bf (4)} $R$ is a Gaussian ring, i.e., $C_{R}(fg) =C_{R}(f)C_{R}(g)$
for any polynomials $f, g$ with coefficients in $R$, where $C_{R}(f)$
is the ideal of $R$ generated by the
coefficients of $f$ called the content ideal of $f$. \\
\indent {\bf (5)} $R$ is a Pr\"ufer ring, i.e., every finitely
generated regular ideal is invertible
(equivalently, every two-generated regular ideal is invertible). \\

\indent In [11], it is proved that each one of the above conditions
implies the following next one (i.e., $(1) \Rightarrow (2)
\Rightarrow (3) \Rightarrow (4) \Rightarrow (5)$), and examples
are given to show that in general the implications can't be
reversed. Moreover, an investigation is carried out to see which
conditions may be added to some of the preceding
properties in order to reverse the implications. \\

\indent Recall that in the domain context, the five classes of
Pr\"ufer-like rings collapse to the notion of Pr\"ufer domain.
From Bazzoni and Glaz [4, Theorem 3.12], we note that a Pr\"ufer
ring $R$ satisfies one of the five conditions if and only if the
total ring of quotients $Tot(R)$ of $R$ satisfies the same
condition. See for instance
[3, 4, 7, 10, 11, 13, 17, 23, 28]. \\

\indent For two rings $A\subseteq B$, we say that $A$ is a module retract (or a subring retract) of $B$ if
there exists an $A$-module homomorphism $\phi: B\longrightarrow A$
such that $\,\phi|_A =id|_A$; $\phi$ is called a module retraction
map. If such a map $\phi$ exists, $B$ contains $A$ as an
$A$-module direct summand. \\
Considerable works have been concerned with the
descent and ascent of a variety of finiteness and related
homological properties between a ring and its subring retract. See for 
instance [5, 6, 9, 18, 19, 24]. \\

\indent A special application of subring retract is the notion of trivial 
ring extension. Let $A$ be a ring,  $E$ an $A$-module and  $R =A \propto  E$, the
set of pairs $(a,e)$ with $a\in A$ and $e\in E$, under
coordinatewise addition and under an adjusted multiplication
defined by $(a,e)(a',e')=(aa',ae'+a'e),$  for all $a,a'\in A,
e,e'\in E$. Then $R$ is called the trivial ring extension of
$A$ by $E$. It is clear that $A$ is a module retract of $R$, where the module
retraction map $\phi$ is defined by $\phi (x,e)=x$.\\
Trivial ring extensions have been studied extensively; the work is
summarized in Glaz [8] and Huckaba  [16]. These extensions have
been useful for solving many open problems and conjectures in both
commutative and non-commutative ring theory. See for instance [8, 16, 20].\\

\indent In this article we investigate the 
transfer of the Pr\"ufer-like properties between a commutative ring and its 
 subring retract. Our results generate new and original 
examples which enrichy the current literature with new families of 
Pr\"ufer-like rings with zerodivisors. \\

\end{section}\bigskip\bigskip

\begin{section}{Pr\"ufer-like properties in subring retract}

 \indent In this section we investigate the transfer
of Gaussian, Pr\"ufer, and arithmetical properties between a ring
and its  subring retract.  \\
 We begin by studying the transfer of Gaussian property. Recall that 
$Nil(R)$ is the set of nilpotent elements in a ring $R$. \\

\begin{thm} 
Let $R$ be a ring and $A$ a subring retract of $R$.   \\
{\bf 1)} If $R$ is a Gaussian ring then so is $A$.  \\
{\bf 2)} Assume that $(A,M)$ is a local ring and $R :=A \propto (A/M)$ be the trivial 
ring extension of $A$ by $A/M$. Then $R$ is a Gaussian ring if and only if 
so is $A$.

\end{thm}

\proof  {\bf 1)} Assume that $R$ is a Gaussian ring and let 
$f(X)=\displaystyle\sum _{i=0}^{n}a_i X^i$, 
$g(X)=\displaystyle\sum _{i=0}^{m}b_i X^i$ be two polynomials of
$A[X]$, where $n$ and $m$ are two positive integers. Our aim is to
prove that $C_{A}(f)C_{A}(g)\subseteq C_{A}(fg)$. Let $\alpha$ be
an element of $C_{A}(f)C_{A}(g)$, we have $\alpha \in
C_{R}(f)C_{R}(g)$ so, since $R$ is Gaussian,  $\alpha \in
C_{R}(fg)$, i.e.,
$\alpha=\displaystyle\sum_{k=0}^{n+m}(\displaystyle\sum_{i+j=k}a_i
b_j)r_k$ where $r_k$ is an element of $ R$ for any $0\leq k \leq
nm$. Then, $\alpha =\phi(\alpha)=
\displaystyle\sum_{k=0}^{n+m}(\displaystyle\sum_{i+j=k}a_i
b_j)\phi(r_k)$ where $\phi$ is the module retraction map, which
prove that $\alpha \in C_{A}(fg)$. Thus, $A$ is a Gaussian
ring. \\

{\bf 2)} Assume that $(A,M)$ is a local ring and $R :=A \propto (A/M)$ be the trivial 
ring extension of $A$ by $A/M$. If $R$ is Gaussian, then so is $A$ by 1). Conversely, 
the fact that $R$ is Gaussian in case $A$ is Gaussian follows easily from the 
characterization of local Gaussian rings given by Tsang ([28]): a local ring $R$ with 
maximal ideal $M$ is Gaussian if and only if for any two elements $a, b$ in $M$ the 
following two conditions hold: 1) $(a,b)^{2} =(a^{2})$ or $(b^{2})$; 
2) if $(a,b)^{2} =(a^{2})$ and $ab =0$, then $b^2 =0$. \cqfd
\bigskip

The necessity of the conditions imposed in Theorem 2.1 will be 
proved in Examples 2.4 and 2.7. \\

Secondly, we study the transfer of Pr\"ufer property between a ring and its subring retract. 
It is clear that each total ring of quotients is a Pr\"ufer ring. Recall that an $R$-module $E$ 
is called a torsion-free if for every regular element $a \in R$ and $e \in E$ such that $ae =0$, 
we have $a =0$ or $e =0$. \\ 

\bigskip

\begin{thm} 
Let $R$ be a ring and $A$ a subring retract of $R$.   \\
{\bf 1)} Assume that the module retraction map $\phi: R\longrightarrow A$  verifies
$Ker(\phi)$ is torsion-free. If $R$ is a Pr\"ufer ring then so is $A$. \\
{\bf 2)} Assume that $(A,M)$ is a local total ring of quotients, where 
$M$ is its maximal ideal; and assume that the module retraction map $\phi$ verifies $Mker(\phi)=0$ 
and $ker(\phi) \subseteq Nil(R)$. 
Then $R$ is a total ring of quotients; in particular, $R$ is Pr\"ufer.
\end{thm}

\proof  {\bf 1)} Assume that $Ker(\phi)$ is torsion-free, where $\phi: R\longrightarrow A$  
is the module retraction map, and $R$ is a Pr\"ufer ring. 
Let $I=\displaystyle\sum_{i=1}^{n}a_i A$ be a finitely
generated regular ideal of $A$ and $a$ be a regular element of
$I$. We tent to prove that $I$ is invertible. Let $b$ be an
element of $R$ such that $ba=0$, we have $\phi(b)a=0$ so
$\phi(b)=0$ since $a$ is a regular element of $A$, i.e. $b\in
Ker(\phi)$. In the other hand, by setting $b =a^{'}+v \in R$ 
where $a^{'} \in A$ and $v \in Ker(\phi)$ (since $Ker(\phi)$ is a direct 
summand of $R$), we obtain that $0 =a^{'}a+va$, and so 
$a^{'}a =va =0$, therefore $a^{'} =0$ and $v =0$ as $a$ is 
regular in $A$ and $Ker(\phi)$ is torsion-free; 
which proves that $a$ is a regular element of $R$ and so the
ideal $J=\displaystyle\sum_{i=1}^{n}a_i R$ is a finitely generated
regular ideal of $R$. Hence, since $R$ is Pr\"ufer, $J$ is
invertible in $R$  and so the polynomial
$f(X)=\displaystyle\sum_{i=1}^{n}a_i X^i$ is Gaussian in $R$
(since $J=C_R (f)$). Using the proof of Theorem 2.1(1), we find that
$f(X)$ is Gaussian in $A$; hence, as $I (=C_{A}(f))$ is a 
regular ideal of $A$, it is invertible
in $A$ (by [3, Theorem 4.2(2)]). Thus, $A$ is Pr\"ufer. \\
{\bf 2)} Assume that $(A,M)$ is a local total ring of quotients, where 
$M$ is its maximal ideal; and assume that the module retraction map $\phi$ verifies $Mker(\phi)=0$ 
and $ker(\phi) \subseteq Nil(R)$. Set $V =Ker(\phi)$. \\
In order to show that $R$ is a total ring of quotients, we have to prove that each element $a+v$ of 
$R$ is invertible or zero-divisor element. Indeed: \\
If $a \in M$, then $a$ is a non invertible element of $A$; 
that's $a$ is zero-divisor in $A$ (since $A$ is a total ring of quotients). Hence 
there exists $b$ nonzero element of $M$ such that
$ab=0$. Therefore, $b(a+v) =0$ as $MV =0$,  which means that $a+v$ is a zero-divisor element in $R$. \\
If $a \notin M$, then $a$ is invertible in $A$ and so in $R$; hence, $a+v$ is invertible in $R$ 
as sum of an invertible element and a nilpotent one. \\
Thus, $R$ is a total ring of quotients.\cqfd
\bigskip



In the following  example we prove that the retraction is not
sufficient to transfer the Pr\"ufer property.\\
\bigskip


\begin{exmp}  Let $(A,M)$ be a non Pr\"ufer local ring and $E$ a nonzero $A$-module such
 that $ME=0$. Let $R :=A\propto E$ be the trivial ring extension of $A$ by $E$. Then: \\
{\bf 1)} $R$ is a total ring of quotients (since $R$ is local with maximal ideal 
$M \propto E$ and $(M \propto E)(0,e) =0$ for each $e \in E$). In particular, $R$ is Pr\"ufer.\\
{\bf 2)} $A$ is a non Pr\"ufer subring retract of $R$.
\end{exmp}
\bigskip

 

In the next example we ensure the necessity of the
conditions imposed in Theorems 2.1(2) and 2.2(2). \\
 \bigskip


\begin{exmp}   Let $(V, M)$ be a rank-one discrete valuation domain such that 
$2\in M$ (for instance, $V:=\Z_{(2)}$). Then  $R:=V \propto~V$ is not Pr\"ufer. 
In particular, $R$ is not Gaussian.
\end{exmp}

\proof  It suffices to show that $R$ is not Pr\"ufer. Let $I :=R(2,0) + R(2,1)$ 
be a finitely generated ideal of $R$. It is clear that $I$ is regular 
(since $(2,0)$ is regular). Since $R$ is local, the $2$-generated regular 
ideal $I$ is invertible if and only if it is principal. Then, again since 
$R$ is local, $I$ is principal if and only if it is generated by one of the two 
generators and this is false, so the conclusion follows easily. \cqfd
\bigskip

We study now the transfer of arithmetical property between a ring
and its  subring retract.  \\. 

\bigskip

\begin{thm}
Let $R$ be a ring and $A$ a subring retract of $R$. If $R$ is an
arithmetical ring then so is $A$.
\end{thm}

\proof By [17, Theorem 2], it suffices to show that 
for any pair of ideals $I$ and $J$ of $A$ such that $I\subseteq J$ and 
$J$ is
finitely generated, there should exist an ideal $H$ of $A$ for which $I=HJ$. \\ 
We have $IR\subseteq JR$
and $JR$ is a finitely generated ideal of $R$; so, as $R$ is
arithmetical, there exists an ideal $L$ of $R$ such that $IR=LJR$
that is $IR=LJ$. Therefore, $I =\phi (IR) =\phi(LJ) =\phi(L)J$ and 
so $A$ is arithmetical.  \cqfd
\bigskip


In the following example we prove that, under the same conditions as 
in Theorem 2.1(2), we can't transfer the arithmetical property from
 $A$ to $R$.  \\
 
 \bigskip
 
\begin{exmp}
 Let $(A,M)$ be a valuation domain which is not a field, where $M$ is its 
maximal ideal. Set $R=A\propto (A/M)$ be the trivial ring extension of $A$ by $A/M$. Then: \\
{\bf 1)} $A$ is an arithmetical subring retract of the local ring $R$. \\ 
{\bf 2)} $R$ is not arithmetical.
\end{exmp}

\proof 1) is clear. Also, we claim that $R$ is not arithmetical. 
Let $I :=R(a,0) + R(0,e)$ be a finitely generated ideal of $R$, where $a$ is 
any nonzero element of $M$ and $e$ is any nonzero element of $A/M$. Since $R$ is local, 
$I$ is principal if and only if it is generated by one of the two generators and 
this is false, so the conclusion follows easily. \cqfd

\bigskip
The following example  proves that the condition ``$A $ is a
subring retract of $R$" can not be removed in the proof of Theorems 2.1(1) and 2.5.\\
\bigskip


\begin{exmp} Let $K$ be a field, $K[X,Y]$ the polynomial ring where $X$ and $Y$ are two indeterminate
elements, and let $Q(K[X])$  be the quotient field of $K[X]$. Then $Q(K[X])[Y]$ is a Pr\"ufer domain 
containing the subring $K[X,Y]$ which isn't a Pr\"ufer domain.   \\
\end{exmp}


\end{section}\bigskip\bigskip

\begin{section}{Applications}

 \indent In this section we give two applications to the results obtained 
in section 2.  The first application is devoted to trivial ring extension 
$R :=A \propto E$ of a ring $A$ by an $A$-module $E$. Recall that 
$A$ is a module retract of $R$, where the module
retraction map $\phi$ is defined by $\phi (x,e)=x$ and $Ker(\phi) =0 \propto E$.\\

\bigskip

\begin{prop} 
Let $A$ be a ring, $E$ an $A$-module and $R :=A \propto E$ be the trivial 
ring extension of $A$ by $E$. Then:   \\
{\bf 1) a)} Assume that $E (=Ker(\phi))$ is torsion-free. If $R$ is a Pr\"ufer ring then so is $A$. \\
{\bf b)} Assume that $(A,M)$ is a local ring, where $M$ is its 
maximal ideal such that $ME =0$.  Then $R$ is a total ring of quotients. In particular, 
$R$ is a Pr\"ufer ring. \\
{\bf 2) a)} If $R$ is Gaussian then so is $A$. \\
{\bf b)} Assume that $E :=A/M$, where $M$ is a maximal ideal of $A$. 
Then $R$ is a Gaussian ring if and only if so is $A$. \\
{\bf 3)} If $R$ is arithmetical then so is $A$. \\
{\bf 4)} $wdim(R) > 1$. \\
\end{prop}

\proof Let's remark first that $R :=A \propto E$, where $(A,M)$ is a local ring and 
$ME =0$, is a total ring of quotients (since $R$ is a local ring with maximal ideal 
$M \propto E$ and $(M \propto E)(0,1) =0_{R}$). By section 2, it remains to show that $wdim(R) > 1$. \\
Let $f \in E-\{0\}$ and $J :=R(0,f) (=0 \propto (Af))$. Consider the exact sequence 
of $R$-modules: \\
$$0\longrightarrow  Ker(u) \longrightarrow  R \buildrel u \over\longrightarrow  J \longrightarrow  0$$
where $u(a,e) =(a,e)(0,f) =(0,af)$. Hence, $Ker(u) =Ann(f) \propto E$. 
We claim that $J$ is not flat. Deny. Then, by [26, Theorem 3.55], 
$J =J \cap Ker(u) =JKer(u) =(0 \propto Af)(Ann(f) \propto E) =0 \propto (Ann(f)f) =0$, 
a contradiction. Therefore, $J$ is not flat and $wdim(R) > 1$. \cqfd

\bigskip

\indent  As shown below, Proposition 3.1 enrichies the literature with new examples of non-Gaussian 
Pr\"ufer rings.   \\

\bigskip


\begin{exmp}  Let $(A,M)$ be a non-Pr\"ufer local domain, where $M$ is 
its maximal ideal, and let $E$ be a nonzero $A$-module such
 that $ME=0$. Let $R :=A\propto E$ be the trivial ring extension of $A$ by $E$. Then: \\
{\bf 1)} $R$ is Pr\"ufer by Proposition 3.1(1.b).  \\
{\bf 2)} $R$ is not Gaussian by Proposition 3.1(2.a) since $A$ is not Gaussian 
(as $A$ is non-Pr\"ufer domain).
\end{exmp}
\bigskip

\indent   For enrich the literature with new examples of non-arithmetical  
Gaussian rings, we propose the next two examples.   \\

\bigskip


\begin{exmp}  Let $(A,M)$ be a valuation domain, where $M$ is 
its maximal ideal, and let $R :=A\propto (A/M)$ be the trivial ring extension of $A$ by $A/M$. Then: \\
{\bf 1)} $R$ is Gaussian by Theorem 3.1(2.b) since $A$ is a valuation domain.  \\
{\bf 2)} $R$ is not arithmetical by Example 2.6(3). \\
\end{exmp}
\bigskip


\begin{exmp}  Let $k$ be a proper subfield of a field $K$ and 
let $R :=k\propto K$ be the trivial ring extension of $k$ by $K$. Then: \\
{\bf 1)} $R$ is Gaussian by [2, Example 2.3(2.b)].  \\
{\bf 2)} $R$ is not arithmetical by [2, Example 2.3(2.c)] since $R$ 
is local. \\
\end{exmp}
\bigskip

\indent  Now we construct an arithmetical ring $R$ such that $wdim(R) > 1$.   \\

\bigskip


\begin{exmp}   Let $K$ be a field and $R :=K\propto K$ be the trivial 
ring extension of $K$ by $K$. Then: \\
{\bf 1)} $R$ is arithmetical.  \\
{\bf 2)} $wdim(R) =\infty $ . \\
\end{exmp}

\proof {\bf 1)} $R$ is arithmetical by [2, Example 2.3(1.a)].  \\
{\bf 2)} The ideal $I :R(0,1)$ 
is not flat by the proof of Proposition 3.1(4). On the other hand, the 
exact sequence of $R$-modules: \\
$$0\longrightarrow  I \longrightarrow  R \buildrel u \over\longrightarrow I \longrightarrow  0$$
where $u(a,e) =(a,e)(0,1) =(0,a)$ shows that $fd_{R}(I) =\infty $. 
Hence, $wdim(R) =\infty $ and this completes the proof. \cqfd
\bigskip

\indent  Let's see in the following example that even if we replace the field 
$K$, in the above example, by a principal total ring of quotients $A$, we don't have in general 
$R :=A \propto A$ is Gaussian; in particular it's not arithmetical. \\
The same example is a Pr\"ufer non Gaussian ring.   \\

\bigskip


\begin{exmp}   Let $A :=\Z/(2^{i}\Z)$, where $i \geq 2$ be an integer, and let 
$R :=A\propto A$ be the trivial ring extension of $A$ by $A$. Then: \\
{\bf 1)} $A$ is a local principal total ring of quotients with maximal ideal $M =2A$.  \\
{\bf 2)} $R$ is a local total ring of quotients. In particular, $R$ is a Pr\"ufer ring.  \\
{\bf 3)} $R$ is not Gaussian. In particular, $R$ is not arithmetical.  \\
\end{exmp}

\proof {\bf 1)} and {\bf 2)} are clear since $R$ is local with maximal ideal 
$M \propto A$ and $(M \propto A)(0,\bar{2^{i-1}}) =0_{R}$. \\
It remains to show that $R$ is not Gaussian. For that let $f :=(\bar{2^{i-1}},0) + (\bar{2^{i-1}},1)X 
\in R[X]$. We have $f^2 =0$ (and so $C_{R}(f^{2}) =0$) and $(C_{R}(f))^2 =R(0,\bar{2^{i-1}}) 
(\not=0_{R}$). Therefore, $R$ is not Gaussian and this completes the proof of Example 3.6. \cqfd

The second application is devoted to the Nagata rings. Let $A$ be a ring and
$R :=A(X)=S^{-1}A[X]$ the localization of $A[X]$ by $S$, where $S$ is
the multiplicative set in $A [X]$ formed by all polynomials 
$f(X)$ such that $C(f) =A$. By construction we have $A(X)=A + XA[X]+ C$ 
where $C=\{ \frac{f(x)}{g(x)}/ f(x),g(x) \in A[X], d^o f(X)<d^o g(x) \;
and\; C(g(x))=A\}$ (by [21, Chapter IV, Proposition 1.4(1)]), 
which implies that $A$ is a module retract of $A(X)$. The ring $R :=A(X)$ is called the Nagata ring. 
See for instance [16, 21]. \\

By section 2, we obtain: \\
\bigskip

\begin{prop} 
Let $A$ be a ring, $R :=A(X)$ be the Nagata ring. Then:   \\
{\bf 1)} Assume that $E (=Ker(\phi))$ is torsion-free. If $R$ is a Pr\"ufer ring then so is $A$. \\
{\bf 2)} If $R :=A(X)$ is Gaussian then so is $A$. \\
{\bf 3)} If $R :=A(X)$ is arithmetical then so is $A$. \\
\end{prop}


\bigskip

\indent  Recall that a subset $S$ of a ring $R$ is called dense if $Ann(S) =0$. A ring 
$R$ is called strongly Pr\"ufer if every finitely generated dense ideal is locally principal. 
Notice that the Nagata ring $A(X)$ is Pr\"ufer if and only if 
$A$ is strongly Pr\"ufer by [16, Theorem 18.10]. For instance, 
 a strongly Pr\"ufer ring is a Pr\"ufer ring by [16, Lemme 18.1]. \\ 
Recall that a ring $R$ satisfies $(CH)$-property if each finitely 
generated ideal of $R$ has a non-zero annihilator. It is clear 
that a $(CH)$-ring is strongly Pr\"ufer ring. For instance, 
the trivial ring extension $R :=A \propto E$ is a $(CH)$-ring (and so 
 strongly Pr\"ufer ring) for each local ring $(A,M)$ (where $M$ is its maximal ideal) and an $A$-module 
$E$ such that $ME =0$.  \\

Now we construct a new examples of non-Gaussian 
Pr\"ufer rings, as shown below.   \\

\bigskip


\begin{exmp}   Let $A$ be a non-Gaussian $(CH)$-ring and let $R :=A(X)$ be the Nagata ring. Then:   \\
{\bf 1)} $R$ is Pr\"ufer.  \\
{\bf 2)} $R$ is a not Gaussian.
\end{exmp}

\proof {\bf 1)} It is clear that $A$ is strongly Pr\"ufer ring 
since $A$ is a $(CH)$-ring. Therefore $R$ is a Pr\"ufer ring 
by [16, Theorem 18.10].  \\
{\bf 2)} $R$ is not Gaussian by Proposition 3.6(2) since $A$ is not Gaussian. \cqfd \bigskip

\noindent {\bf ACKNOWLEDGEMENTS.} The authors would like to express their sincere thanks for 
the referee for his/her helpful suggestions and comments, which have greatly improved this paper. \\

\end{section}\bigskip\bigskip



\bigskip\bigskip

\end{document}